\documentclass[12pt,reqno,oneside]{amsart}
\usepackage{amssymb}

\newtheorem{theorem}{Theorem}
\newtheorem{proposition}[theorem]{Proposition}
\newtheorem{lemma}[theorem]{Lemma}

\newcommand\lp{\left(}
\newcommand\rp{\right)}

\newcommand\pp{p}
\newcommand{\hu}{\hat{u}}
\newcommand{\hv}{\hat{v}}

\newcommand\Res{\mathop{\mathrm{Res}}}

\newcommand\esp{e}

\newcommand\natnums{\mathbb{N}}

\newcommand\thxbod{The author extends his thanks to Prof. D.  Jackson
    for his helpful suggestions during the preparation of this
    article.}

\title{Composition sums related to the hypergeometric function}
\author{R. Milson}
\address{Dept of Math., Dalhousie U., Halifax, Canada, B3J 3J5}

\curraddr{McGill University, Montreal}
\email{milson@math.mcgill.ca}
\thanks{\thxbod}

\begin{document}
\maketitle
%
%
%
\begin{abstract}
  The present note considers a certain family of sums indexed by the
  set of fixed length compositions of a given number.  The sums in
  question cannot be realized as weighted compositions.  However they
  can be be related to the hypergeometric function, thereby allowing
  one to factorize the corresponding generating polynomials.  This
  factorization leads to some interesting identities.
\end{abstract}
\bigskip

\pagestyle{myheadings}
\markboth{}{Composition sum identities}

An $l$-composition of a natural number $n$ is an ordered list of $l$
positive integers $\pp=(p_1, p_2, \ldots, p_l)$ such that $p_1+\ldots
+ p_l=n$.  The purpose of the present note is to exhibit closed form
expressions for a certain family of sums indexed by the set of fixed
length compositions of a given number.  There are examples of such
identities relating composition sums to Stirling numbers
\cite{Sitgreaves,HomenkoStrok} and to Fibonacci numbers
\cite{HoggattLind, MoserWhitney}.  It should be noted that the sums
introduced in the preceding references can all be regarded as
enumerations of weighted compositions \cite{MoserWhitney}.  To be more
specific, a weighted composition is one where the $j^{\text{th}}$ term
is given a weight $w_j$, and the enumeration of such compositions is
defined as
$$S(l,n)=\sum w_{p_{1}} w_{p_{2}} \ldots w_{p_l},$$
where the sum is
taken over all $l$-compositions of a fixed $n$.  It therefore follows
that all such sums can be realized in terms of a certain type of
generating function:
$$\frac{1}{1-tw(x)}= \sum_{l,n} S(l,n)t^l x^n,$$
where $w(x)=w_1 x+w_2
x^2+ \ldots$.

The same cannot be said of the sums introduced in the present note.
Instead, the identities to be discussed here come about because one
can relate the sums in question to the hypergeometric function,
$F(\alpha,\beta,\gamma; z)$.  This is accomplished by
gauge-transforming a certain second-order differential equation into
the hypergeometric equation , and then taking the residue with respect
to the $\gamma$ parameter.  For a number of other combinatorial
identities involving hypergeometric functions, as well as for other
properties of these fascinating mathematical objects the reader is
referred to \cite{Bateman}.

For a given composition, $\pp$, let $L(\pp)$
and $R(\pp)$ denote, respectively, the products of the left and right
partial sums of $\pp$.  To wit
\begin{align*}
  L(\pp)& = p_1 (p_1+p_2) \ldots (p_1+\ldots + p_{l-1})\, n,\\
  R(\pp)& = p_l (p_{l-1}+p_l) \ldots (p_2+\ldots + p_l)\, n.
\end{align*}
Throughout let $\esp_j(x_1,\ldots,x_k)$ denote the $j^{\text{th}}$
elementary symmetric function of the $x$'s, i.e. the coefficient of
$X^{k-j}$ in the expansion of $(X+x_1)\ldots(X+x_k)$.

The weighted sums that will be discussed here are indexed by three
natural numbers, and are defined by
$$S(k,l,n) = \sum \frac{\esp_k(p_1,\ldots,p_l)}{L(\pp)R(\pp)},$$
where the sum is taken over all $l$-compositions of $n$.  For every
$n>0$ define a corresponding generating polynomial by
$$P_n(u,v) = \sum_{k=0}^l\sum_{l=1}^n S(k,l,n) u^k v^{l-k}.$$
The main
result of the present note is the following factorization of this
generating polynomial.
\begin{theorem} 
  In the even case, say $n=2m$, one has
  $$(n!)^2 P_n(u,v)= \prod_{i=0}^{m-1} \left[ (u+v+q_i)(u+v+q_{i+1}) +
    r_i u\right],
  $$
  where $q_i=i(n-i)$ and $r_i = (n-1-2i)^2$.  In the odd case, say
  $n=2m+1$ one has
  $$(n!)^2 P_n(u,v)=(u+v+q_m) \prod_{i=0}^{m-1} \left[
    (u+v+q_i)(u+v+q_{i+1}) + r_i u\right].
  $$
\end{theorem}

The proof of the theorem will be given below.  First it is worth
remarking that the above theorem implies some attractive identities:
\begin{equation}
\label{eqn:id1}
\sum \frac{p_1 p_2 \ldots p_l}{L(\pp) R(\pp)} 
= \frac{1}{n}\,\esp_{l-1}\! \lp \frac{1}{1\cdot 2} ,
\frac{1}{2\cdot 3}
, \ldots , \frac{1}{(n-1)\cdot n}\rp,
\end{equation}
\begin{multline}
\label{eqn:id2}
\sum \frac{(p_1-1) (p_2-1) \ldots (p_l-1)}{L(\pp) R(\pp)} \\
= \frac{n-1}{n^2}\,\esp_{l-1}\! \lp \frac{r_1}{q_1 q_2},
\frac{r_2}{q_2 q_3} ,\ldots, \frac{r_{m-1}}{q_{m-1} q_m}\rp,
\end{multline}
where the sums are taken over all $l$-compositions of $n$, and where
$m$ in equation \eqref{eqn:id2} is the largest integer smaller or
equal to $n/2$.

To prove the identity in \eqref{eqn:id1} note that the left hand side
is just $S(l,l,n)$, and hence can be obtained as the coefficient of
$u^l$ in the factorizations of Theorem 1.  An easy calculation shows that
$$q_i+ q_{i+1}+ r_i = i(i+1)+(n-i-1)(n-i).$$
Hence 
$$P_n(u,0) = \frac{u}{n} \prod_{i=1}^{n-1} \lp \frac{u}{i(i+1)}+1\rp$$
and the desired identity follows immediately.

The identity in \eqref{eqn:id2} is proved by noting that the left hand
side is given by the alternating sum
$$
\sum_{k=0}^l (-1)^{l-k} S(k,l,n),$$
and hence can be obtained as
the coefficient of $u^l$ in $P_n(u,-u)$.  Using Theorem 1 to evaluate
the latter leads to \eqref{eqn:id2}.

Using Theorem 1 to evaluate $P_n(0,v)$ one also obtains the identity
$$
\sum \frac{1}{L(p)R(p)} = \frac{1}{n^2}\,
\esp_{l-1}\!\lp\frac{1}{q_1},\frac{1}{q_2},\ldots,
\frac{1}{q_{n-1}}\rp.$$
However, this identity is considerably less
interesting, because it can be obtained by a simple rearrangement of
factors in the summands of the left hand side in question.

The proof of Theorem 1 will require a number of intermediate results.
Let $a_{ij},\, i,j\in\natnums$ be indeterminates.  For a composition
$p_1,p_2,\ldots,p_l$ of $n$ let $s_j$ denote the $j^{\text{th}}$ left
partial sum, $p_1+\ldots+p_j$, and set $a(p)=a_{0 s_1} a_{s_1 s_2}
a_{s_2 s_3} \ldots a_{s_{l-1} n}$.
\begin{lemma}
 \label{lemma:ordpart}
 Defining $f_n,\, n\in\natnums$ recursively by
 $$f_n = \sum_{j=0}^{n-1} a_{jn} f_j,\quad f_0=1,$$
 one has $f_n =
 \sum_p a(p)$ where the sum is taken over all compositions (of all
 lengths) of $n$.
\end{lemma}

For $f(\gamma)$, a rational function of an indeterminate $\gamma$,
let $\Res(f;\gamma=\gamma_0)$ denote the residue of this function at
the value $\gamma_0$, i.e. the coefficient of $(\gamma-\gamma_0)^{-1}$
in the Laurent series expansion of $f$ about $\gamma=\gamma_0$.  Let
$$f(u,v,\gamma;z)=1+\sum_{n>0} f_n(u,v,\gamma)z^n$$
denote the unique
formal power series solution of
\begin{equation}
  \label{eq:fdef}
z^2 f_{zz} + (1-\gamma) z f_z + \lp \frac{vz}{1-z} +
\frac{uz}{(1-z)^2}\rp f = 0,\quad f(0)=1. 
\end{equation}
\begin{proposition}
\label{prop:res}
  $n P_n(u,v) = \Res(f_n; \gamma=n)$.
\end{proposition}
\begin{proof}
  Rewriting $z\,(1-z)^{-1}$ as $\sum_{n>0} z^n$ and $z\,(1-z)^{-2}$ as
  $\sum_{n>0} nz^n$ one obtains the following recurrence relation for
  the coefficients of $f$:
  $$n(\gamma-n) f_n = \sum_{i=0}^{n-1} ((n-i)u+v)f_i.$$
  Set
  $$a_{ij} = \frac{(j-i)u+v}{j(\gamma-j)}$$
  and note that for a
  composition $p_1,\ldots,p_l$ of $n$ one has
  \begin{align*}
    a(p) &= \frac{p_1 u+v}{s_1(\gamma-s_1)}\times
    \frac{p_2u+v}{s_2(\gamma-s_2)} \times \ldots\times
    \frac{p_{l-1}u+v}{s_{l-1}(\gamma-s_{l-1})}\times
    \frac{p_lu+v}{n(\gamma-n)} \\
    &= \frac{\sum_{k=0}^l \esp_k(p_1,\ldots,p_l) u^k v^{l-k}} {L(p)
      (\gamma-s_1)(\gamma-s_2)\ldots(\gamma-n)}.
  \end{align*}
  Taking the residue of the right hand side at $\gamma=n$ and applying
  Lemma \ref{lemma:ordpart} gives the desired conclusion.
\end{proof}
Let $\alpha, \beta, \gamma$ be indeterminates, and let
$F(\alpha,\beta,\gamma;z)$ denote the usual hypergeometric series
$$\sum_{n=0}^\infty \frac{(\alpha)_n\,(\beta)_n}{n!\,(\gamma)_n}\,z^n
,$$
where $(x)_n$ is an abbreviation for the expression
$x(x+1)\ldots(x+n-1)$.
\begin{proposition} 
\label{prop:hypg}
Setting
\begin{align}
  \hu&=\frac14\,(\alpha+\beta+\gamma)(2-\alpha-\beta-\gamma),\\
  \nonumber \hv&=\frac14\,(\alpha-\beta-\gamma)(\alpha-\beta+\gamma).
\end{align}
one has
\begin{equation}
  \label{eq:fFrel}
f(\hu,\hv,\gamma;z)=
(1-z)^{(\alpha+\beta+\gamma)/2}\,F(\alpha,\beta,1-\gamma;z).
\end{equation}
\end{proposition}
\begin{proof}
  Substituting \eqref{eq:fFrel} into \eqref{eq:fdef} yields the
  following equation for $F$
  $$
  z^2 F_{zz} + (1-\gamma)zF_z - (\alpha+\beta+\gamma) F_z - \alpha\beta
  F = 0.$$
  Multiplying the above by $(1-z)/z$ yields the usual
  hypergeometric equation with parameters $\alpha, \beta, 1-\gamma$.
\end{proof}
\begin{proof}[Proof of Theorem 1.]
  Expanding $(1-z)^{(\alpha+\beta+\gamma)/2}$ into a power series in
  $z$, and using Proposition \ref{prop:hypg} one obtains that
  $$f_n(\hu,\hv,\gamma) =
  \frac{(\alpha)_n\,(\beta)_n}{n!\,(1-\gamma)_n}+\ldots,$$
  where the
  remainder is a finite sum of rational functions in $\alpha, \beta,
  \gamma$ which do not contain a factor of $\gamma-n$ in the
  denominator.  Hence, by Proposition \ref{prop:res}
  $$P_n(\hu,\hv) = (-1)^n\,\frac{(\alpha)_n\,(\beta)_n}{(n!)^2}.$$
  An
  elementary calculation shows that for $0\leq i<n/2$
  $$
  (\hu+\hv+q_i)(\hu+\hv+q_{i+1})\, + \,r_i\hu =
  (\alpha+i) (\beta+i) (\alpha+n-1-i)(\beta+n-1-i)
 $$   
 and that for $n=2m+1$
 $$
 \hu+\hv+q_m = -(\alpha+m)(\beta+m).$$
 The desired factorizations
 follow immediately.
\end{proof}

\end{document}